\documentclass[a4paper,11pt]{amsart}
\usepackage[english]{babel} 
\usepackage[leqno]{amsmath}
\usepackage{amssymb,amsthm}
\usepackage{amscd}
\usepackage{enumerate}

\newtheorem{thm}{Theorem}[subsection]
\newtheorem{lemma}[thm]{Lemma}

\newtheorem{prop}[thm]{Proposition}

\theoremstyle{remark}
\newtheorem{remark}[thm]{Remark}
\theoremstyle{definition}
\newtheorem{defi}[thm]{Definition}
\newtheorem{nota}[thm]{}

\newtheorem{example}[thm]{Example}

\newcommand{\la}{\longrightarrow}
\newcommand{\ha}{\hookrightarrow}
\newcommand{\da}{\dashrightarrow}

\newcommand{\ov}{\overline}

\newcommand{\Spec}{\operatorname{Spec}}

\newcommand{\im}{\operatorname{Im}}

\newcommand{\Pic}{\operatorname{Pic}}
\newcommand{\supp}{\operatorname{Supp}}
\newcommand{\mdeg}{\operatorname{{\underline{de}g}}}

\def\E{\mathcal E}
\def\L{\mathcal L}

\def\O{\mathcal O}

\def\X{\mathcal X}

\newcommand{\Z}{\mathbb{Z}}

\def\NN{\operatorname{N}}

\def\dcg{\Delta _X}

\newcommand{\G}{\Gamma }

\def\md{\underline{d}}

\newcommand{\g}{\gamma}

\newcommand{\dd}{\delta}

\newcommand{\eSd}{\epsilon^{\underline d}_S}

\newcommand{\pr}[1]{\mathbb{P}^{#1}}
\newcommand{\Xn}{X^{\nu}}

\newcommand{\mgbar}{\ov{M}_g}

\def\ad{\alpha^d_f}
\def\adX{\alpha^d_{f,X}}
\def\agX{\alpha^{g-1}_{f,X}}
\def\ag{\alpha^{g-1}_f}

\def\Lg{{\mathcal L}_{g-1}}

\newcommand{\dX}{\dot{\X}}

\def\nf{N^d_f}

\def\PXb{\overline{P^d_X}}
\def\vine{X_{\delta}}
\def\PXbd{\overline{P^d_{X_{\delta}}}}

\def\PSd{P^{\md}_S}

\def\pXgb{\overline{P^{g-1}_X}}
\def\PXgb{\overline{P^{g-1}_X}}
\def\pf{P^d_f}

\def\pfb{\overline{P^d_f}}
\def\pfgb{\overline{P^{g-1}_f}}

\def\Wmd{W_{\md}(X)}

\def\Amd{A_{\md}(X)}
\def\Xmd{X^{\md}}
\def\Xmds{X^{\md}_{\sigma}}

\def\amd{\alpha^{\md}_X}
\def\amds{\alpha^{\md}_{X,\sigma}}

\def\BXs{\Sigma^{ss} (X)}
\def\BX{\Sigma(X)}

\def\BXS{\Sigma (\XS)}

\newcommand{\picf}[1]{\Pic_f^{#1}}
\newcommand{\picX}[1]{\Pic^{#1}X}

\def\pdbst{\overline{\mathcal{P}}_{d,g}}
\def\pdbstn{\pdbst^{Ner}}
\def\pdbn{\overline{P}_{d,g}^{Ner}}

\def\mgbst{\overline{\mathcal{M}}_g}

\newcommand{\tw}{\operatorname{Tw}}
\newcommand{\twf}{\tw_f}

\newcommand{\XS}{X_S^{\nu}}

\newcommand{\hX}{{\widehat{X}}}
\newcommand{\hXS}{{\widehat{X}}_S}

\newcommand{\gen}{\X_K}

\newcommand{\sing}{X_{\text{sing}}}

\begin{document}

\centerline{\bf{Compactified Jacobians, Abel maps and Theta divisors}}

\centerline{Lucia Caporaso}

\


\tableofcontents
\

\section{Introduction}

\noindent{\it Outline of the paper.}
This paper is an expository account about compactified Picard schemes of nodal curves
 and some related topics.
After some preliminaries,
 N\'eron models are used to classify different compactified Picard schemes,  
Abel maps are studied accordingly, and finally  some recent results on the Theta divisor
 are reviewed.

Several examples are included, as an attempt to elucidate
some important  parts that   require a consistent amount of
technical work to be  rigorously settled. 
To fill out the  unavoidable gaps,
various references are given throughout the paper.

\

\noindent{\it Aknowledgements.}
A large part of this paper is based on a talk that I gave at the conference on Curves and Abelian varieties,
held in Spring 2006  to
celebrate  Roy Smith's   birthday. I wish to dedicate this paper to him.
I also want to warmly thank  Valery Alexeev and Elham Izadi for organizing the conference,
and Eduardo Esteves for carefully reading a preliminary version.

\

\noindent{\it Conventions.}
We fix  an algebraically closed field $k$ and  work with schemes  locally of finite type over $k$, unless
differently specified.
 
$X$ will always be a connected, reduced, projective  curve over $k$, having at most nodes as singularities.
We denote by $g$ the arithmetic genus of $X$, by $\dd$ the number of its nodes and by 
$\gamma$ the number of its
irreducible components. The dual graph of $X$ (having as vertices  the $\gamma$
irreducible components of $X$ and as edges the $\dd$ nodes of $X$) is denoted $\Gamma_X$.

$X$ is called {\it of compact type} if $\Gamma_X$ is a tree; more generally, $X$ is called {\it tree-like}
if   $\Gamma_X$  becomes a tree after  all loops are removed (in particular every 
irreducible curve is tree-like).

By $f:\X \to B$ we denote a  {\it one-parameter smoothing} of $X$. That is,
$B=\Spec R$ is the spectrum of a discrete valuation ring $R$ having residue field $k$ and quotient field $K$;
the closed fiber of $f$ is $X$ and the generic fiber, denoted $\gen$, is a smooth projective curve over $K$.
Everything we shall say holds (mutatis mutandis) if one replaces $B$ by any Dedekind scheme.

The total space $\X$ is a normal surface with singularities of type $A_n$ at the nodes of $X\subset \X$.
If $\X$ is nonsingular, we will say that $f$ is a {\it regular  smoothing}.

\section{Compactified Picard schemes}
\subsection{Set up}

\begin{nota}{\it Jacobians, their torsors and their models.}
We consider the Picard scheme of $\gen$
$$
\Pic\gen=\coprod _{d\in \Z}\Pic^d\gen
$$
where $\Pic^d\gen$ is the variety parametrizing isomorphism classes of line bundles of degree $d$ on $\gen$.
In particular $\Pic^0\gen$ is an abelian variety over $K$ and $\Pic^d\gen$ a torsor under it.

What about   models  of $\Pic^d\gen$ over $B$?  The analysis of   such models 
 (by which   we mean  integral schemes
over $B$ whose generic fiber is $\Pic^d\gen$)
is
a key to understand the properties of any {\it compactified Jacobian} or 
{\it compactified Picard scheme}  of  $X$, which is one of the main
themes of this paper. 

We are going to introduce three different types of models for $\Pic^d\gen$ over $B$. 
The first two
(the Picard scheme and the N\'eron model) are  uniquely determined by their defining properties, the third  
 consists of a miscellany of different models (compactified Picard schemes, see \ref{comppic}).
We shall see how they relate to one another; in particular,
they are all isomorphic if and only if $X$ is nonsingular.
\end{nota}

\begin{nota}{\it The (relative) Picard scheme.}
\label{relpic}
For every   $f:\X \to B$ and for every $d$  there exists the
 (relative, degree $d$) Picard scheme    
$$
\lambda_d:\picf{d}\to B
$$  
(often denoted $\Pic^d_{\X/B}\to B$). The existence and basic theory are due to 
A. Grothendieck, P. Deligne and D. Mumford (see \cite{deligne}, \cite{SGA}, \cite{mumford});
we  refer to \cite{BLR} for a unified account.
Over every point in $B$, the fiber of $\lambda_d$   is the variety of isomorphism classes of line bundles
of degree $d$ on the corresponding fiber of $f$. So, the generic fiber of $\lambda_d$ is $\Pic^d\gen$ and the
closed fiber is $\Pic^dX$  (see (\ref{pX}) for an explicit description of $\Pic^dX$).

The moduli property of the Picard scheme is 
expressed in Proposition 4, section 8.1. p. 204 of \cite{BLR}; loosely speaking  it amounts to
the following. For every $B$-scheme   $T\to B$,
denote $f_T:\X_T:=\X\times _B T\to T$ the base change of $f$. For every line bundle $\L$ of relative degree
$d$ on
$\X_T$ there exists a unique ``moduli morphism"  $\mu_{\L}:T\to \picf{d}$ mapping $t\in T$ to the isomorphism class
of the restriction of $\L$ to $f_T^{-1}(t)$;
for any $M\in \Pic T$, we have    $\mu_{\L}=\mu_{\L\otimes f_T^*M}$. Conversely, given a morphism
$\mu:T\to \picf{d}$,
the obstruction to the existence of a line bundle $\L$ on $\X_T$ having $\mu$ as moduli map  lies
in the Brauer group of
$T$. If
$f$ has a section,   there is no   obstruction.

The union $\Pic _f:=\coprod_{d\in \Z}\picf{d} $ is a group scheme over $B$ with respect to tensor product.

An important fact is that {\it $\Pic_f^d\to B$ is separated if and only if $X$ is irreducible} (see \ref{nosep}).
\end{nota}

\begin{nota}{\it The generalized jacobian.}
\label{genjac}
The generalized jacobian of $X$,   denoted $J(X)$, parametrizes isomorphism classes of line
bundles having degree $0$ on every irreducible component of $X$; thus $J(X)$ is a commutative algebraic group with
respect to tensor product.
It is well known that $J(X)$  is a semiabelian variety, i.e. there exists an exact sequence
\begin{equation}
\label{gj}0\la (k^*)^b\la  J(X)\stackrel{\nu^*}{\la}  J(\Xn) \la 0
\end{equation}
where $\nu:\Xn\to X$ is the  normalization of $X$ 
(hence $J(\Xn)$ is an abelian variety) and 
$$
b=\dd-\g+1=b_1(\Gamma _X)
$$ is the first betti number of the dual graph of $X$,
$\Gamma _X$. As a
consequence we see that {\it $J(X)$ is projective if and only if $X$ is a curve of compact type}.

To relate the generalized jacobian to the Picard scheme, note that for every $d\in \Z$ we have
\begin{equation}
\label{pX}
\Pic^dX=\coprod_{\md\in \Z^{\g}}\picX{\md}
\end{equation}
where   $\picX{\md}$ parametrizes line bundles of multidegree $\md$ on $X$. For example
$\Pic^{(0,\ldots,0)}X=J(X)$. Of course $\picX{\md}$ is non-canonically
isomorphic to $J(X)$ for every $\md$.
We shall sometimes abuse terminology by calling $\picX{\md}$ a generalized jacobian.

For every   $f:\X \to B$ there exists the relative jacobian,
$ 
J_f\to B
$ 
(often denoted $J_{\X/B}\to B$) which is a group scheme over $B$ having fibers  the generalized jacobians of the
fibers of $f$.

The connected component of the identity, $(\Pic_f)^0\to B$, of the group scheme $\Pic_f\to B$ is canonically
identified with 
$ J_f\to B.$ 
\end{nota}
\begin{nota}{\it The N\'eron model}
\label{neron} We here describe some well known results of A. N\'eron and M. Raynaud
(\cite{neron}, \cite{raynaud}) on the existence of N\'eron models and on N\'eron models of Picard varieties.
We  refer to \cite{BLR} for all details, and to \cite{artin} for a    synthesis
of the basic theory.
The N\'eron model of $\Pic^d\gen$  will here be denoted  
$$
\sigma_d:\nf \la B
$$
(a common notation is  $\NN(\Pic^d\gen)\to B$); $\sigma_d$ is a  smooth, separated morphism of finite type.
If $d=0$ then $N^0_f\to B$ is a group scheme whose identity component is  $ J_f\to B$;
for a general $d$,
$\nf \to B$
is a torsor under $N^0_f\to B$
(the N\'eron model of $\Pic^0\gen$).
The N\'eron model is uniquely determined by  the {\it N\'eron mapping property}  (\cite{BLR}, Def.~1, p. 12), which
is the following. For every   scheme $Z$ smooth over $B$, every map $u_K$ from its generic fiber  
to $\Pic^d\gen$
 extends uniquely to a regular $B$-morphism  
$N(u_K):Z\to N^d_f$. 

The closed fiber of $N^d_f$,  viewed simply as a scheme (forgetting its torsor structure),
depends only on the type of singularities of the surface $\X$, in particular it does not depend on $d$; so,
if $f:\X \to B$ is a regular smoothing, we will denote its
closed fiber
$N_X^d$, or $N_X$ when no confusion is possible. 

The scheme $N_X$ can be described in various ways; we begin using combinatorics, 
as in \cite{OS} sections 4 and
14, to which we refer for more details. Consider the dual graph $\Gamma_X$ of $X$, let $c(X)$ be the {\it
complexity} of
$\G_X$; recall that the complexity of a graph   is the number of 
trees contained in it and passing through all of its vertices (i.e. the number of so-called ``spanning trees"). 
Now,
$N_X$ is the disjoint union of
$c(X)$ copies of the generalized jacobian of $X$.
Therefore {\it $N_X$ is irreducible (and isomorphic to $J(X)$) if and only if  $X$ is
a tree-like curve.}
An alternative description  of $N_X$ will be given in the sequel.
\end{nota}
\begin{example}{\it The curve $\vine$.}
\label{vine} 
Denote by 
$$
\vine=C_1\cup C_2,\  \  \  \  \delta=\#(C_1\cap C_2)=\#\sing
$$
the union of two smooth curves meeting transversally at $\dd$ points.   
We call $g_i$ the genus of $C_i$ so that the genus $g$ of $\vine$ is $g=g_1+g_2+\dd -1$.
$\vine$  is sometimes called a {\it vine} curve;
it will be our leading example throughout the paper.

$\Gamma_{\vine}$ consists of two vertices joined by $\dd$ edges; therefore every edge is a spanning tree and 
$
c(X)=\dd 
$. We obtain that $N_X$ is the disjoint union of $\dd$ copies of $J(X)$.
\end{example}

\subsection{Relation between Picard scheme and N\'eron model}
\label{compNP}
\begin{nota}{\it A canonical quotient.}
\label{}
Fix $f:\X \to B$ a regular (for simplicity) smoothing of $X$ and consider the two models
of $\Pic^d\gen$ that we have introduced so far:
the Picard scheme  $\lambda_d:\Pic^d_f \to B$ and the N\'eron model $\sigma_d:\nf \to B$.
The N\'eron mapping property yields a canonical regular $B$-map
\begin{equation}
\label{qf}
q_f:\Pic^d_f\to \nf
\end{equation}
 extending the identity on the generic fibers
(if $d=0$ our $q_f$ is the map ``Ner" in diagram 1.21 of \cite{artin}); we omit    $d$ for simplicity. 
$q_f$ is a surjection, indeed if $d=0$ it is a quotient of group schemes.
So $\nf$ is sometimes called the ``largest separated quotient of the degree-$d$ Picard scheme".
The
restriction  of $q_f$   to the closed fibers depends on $f$; 
we shall now concentrate on it.
\end{nota}
\begin{nota}{\it Twisters.}
\label{tw}
To a regular smoothing $f$ of $X$ we can  associate   a discrete subgroup $\twf X$ of $\Pic^0X$;
$\twf X$ is the set of all line bundles of the form $\O_{\X}(D)_{|X}$ where $D$ is a divisor on $\X$ supported on  
the closed fiber $X$. Elements of $\twf X$ are called {\it twisters} (or $f$-twisters).
The multidegree map 
$$
\mdeg :\twf X \la \Z^{\gamma}
$$
has image a group called $\Lambda_X$:
$$
\Lambda_X=\mdeg (\twf X) \subset \{\md \in \Z^{\gamma} : |\md|=0\}\subset \Z^{\gamma}.
$$
\end{nota}
\begin{remark}
\label{dep}
Observe that while  $\twf X$ depends on $f$ (unless $X$ is   tree-like)
$\Lambda_X$ does not. For example, if $X=C_1\cup C_2$ with $\#C_1\cap C_2=\dd \geq 2$, for every $n\neq 0$
we have that there exist regular smoothings $f:\X \to B$ and $f':\X' \to B$
of $X$ such that $\O_{\X}(nC_1)_{|X}\not\cong\O_{\X'}(nC_1)_{|X}$.
On the other hand, for any $f$ we have $\mdeg\O_{\X}(nC_1)_{|X}=(-n\dd,n\dd)$.
\end{remark}

\begin{nota}{\it Towards separatedness.}
\label{nosep}
Consider for a moment the case $d=0$; recall that $N^0_f$ is a separated model of $\Pic^0\gen$ endowed with a
universal mapping property. 
Consider   the line bundles $\O_{\X}$ and $\L:=\O_{\X}(D)$ 
for any divisor $D$ on $\X$ such that $\supp D\subset X$.
Suppose   that $D$ is not a multiple of $X$,
which is equivalent to ask that $\mdeg \O_{\X}(D)_{|X}\neq (0,\ldots,0)$. 
Thus $\L$ and $\O_{\X}$ determine two different sections,  
 $\mu_{\O_{\X}}:B\to \Pic^0_f$ and $\mu_{\L}:B\to \Pic^0_f$,
of the Picard scheme $\Pic^0_f\to B$ (cf. \ref{relpic}). These two sections 
of course coincide on the generic point $\Spec K$.

The generic fiber of $N^0_f$ is the same as that of $\Pic^0_f$, so we may interpret
$(\mu_{\O_{\X}})_{|\Spec K}=(\mu_{\L})_{|\Spec K}$ as a map from $\Spec K$ to the generic fiber  of $N^0_f$. 
By the N\'eron mapping property, there exists a unique morphism 
$\ov{\mu}:B\to N^0_f$ extending   $(\mu_{\O_{\X}})_{|\Spec K}$ and $(\mu_{\L})_{|\Spec K}$
and such that 
$$\ov{\mu}=q_f\circ  \mu_{\O_{\X}}=q_f\circ \mu_{\L}.$$ 
This implies that  in the closed fiber of $N^0_f$ there is a  unique point corresponding to 
both $(\O_{\X})_{|X}=\O_X$ and $\O_{\X}(D) _{|X}$ for any $D$; in other words   $q_f$ maps    $\twf X$
to one point.
More generally, we have 
\begin{lemma}
\label{} Let $L, L'\in \Pic^dX$ and assume that there exists a  smoothing  $f:\X\to B$ of $X$ and a $T\in \twf
X$
such that $L'\cong L\otimes T$. Then $q_f(L)=q_f(L')$.
\end{lemma}
\begin{proof} (The lemma and its proof hold for every smoothing $f$ of $X$, regardless of it being  regular). As
$T\in
\twf X$ there exists a line bundle $\mathcal T=\O_{\X}(D)$ on $\X$ which restricts to $T$ on the closed fiber $X$. 
Suppose that there exists   $\L\in \Pic\X$ restricting to $L$ on $X$; set $\L':=\L \otimes {\mathcal T}$
so that $\L'_{|X}=L'$ and $\L_{|\gen}=\L'_{|\gen}$.
We
can argue as we did in \ref{nosep}   (with respect to the pair $\L$ and $\O_{\X}$) 
with respect to the pair $\L$ and $\L '$. So we are done. 

Now, up to replacing $f$ with an \'etale base change, we can assume that such a line bundle $\L$ exists.
Since the formation of N\'eron models commutes with \'etale base change, we are done.
\end{proof}
\end{nota}
\begin{nota}{\it Multidegree classes.}
\label{dcg}
The previous  discussion motivates the following definition.
Let $\md$ and $\md'$ be two multidegrees (so that $\md, \md'\in \Z^{\gamma}$); we define them to be equivalent
if their difference is the multidegree of a twister, i.e.:
$
\md \equiv \md ' \  \Leftrightarrow \  \md - \md '\in \Lambda_X
$.
Now the quotient of the set of   multidegrees with fixed total degree by this equivalence relation is a
finite set $\dcg^d$:
$$
\dcg^d:=\frac{\{\md \in \Z^{\gamma}: |\md|=d\}}{\equiv}.
$$
It is well known that the cardinality of
$\dcg^d$ is independent of $d$
(so that we shall sometimes simply write  $\dcg$) and it is equal to $c(X)$
(defined in \ref{neron}).
$\dcg^d$  naturally
labels the connected/irreducible components of $N_X^d$; indeed we have
\begin{equation}
\label{}
N_X^d=\coprod_{\mu\in \dcg^d}N_{X,\mu}
\end{equation}
with non canonical isomorphisms
$$
N_{X,\mu}\cong J(X),\  \  \forall \mu\in \dcg^d.
$$

Finally, the restriction of $q_f$ to the closed fiber $\Pic^dX$ of $\lambda_d$
is surjective and induces an isomorphism of each connected component $\Pic^{\md}X$ 
of $\Pic^dX$ with the connected
component of $N_X^d$ corresponding to the class $\mu$ of $\md$ in $\dcg^d$:
$$
(q_f)_{|\picX{\md}}:\picX{\md}\stackrel{\cong}{\la}N_{X,\mu},\  \  \text{where } \mu=[\md]\in\dcg^d
$$
(more details in  chapter 9 of \cite{BLR} or in  \cite{artin}).
\begin{example}{\it Structure of $\Delta_{\vine}^d$.}
\label{dcgvine}
In the  situation of Example~\ref{vine},
one has that the group of multidegrees of twisters is
$$
\Lambda_{\vine}=\{(n\dd,-n\dd), \forall n\in \Z\}\subset \Z^2
$$
(see remark~\ref{dep}). Hence $\Delta_{\vine}^d$ has cardinality $\delta$ (indeed $\Delta_{\vine}^0\cong \Z/\delta
\Z$).
\end{example}
\end{nota}
\subsection{Types of Compactified  Picard schemes}

\begin{nota}{\it A general notion of compactified Picard scheme.}
\label{comppic}
In this paper, a  (degree-$d$) compactified Picard scheme  for $X$  is a projective, reduced
scheme 
$\PXb$ such that
for every $f:\X \to B$    
there exists a unique integral scheme $\pfb$ with a projective morphism 
\begin{equation}
\label{pd}
\pi_d:\pfb\to B\end{equation}
  whose generic fiber is $\Pic^d\gen$ and whose closed fiber
is $\PXb$. We call $\pfb$ the (relative, degree-$d$) compactified Picard scheme associated to $f$.
We denote by $\pf\subset \pfb$ the smooth  locus of $\pi_d$.

A compactified Picard scheme is also endowed with some geometric meaning. More precisely,
$\PXb$ and $\pfb$ will be (coarse or fine) moduli schemes for certain functors
strictly related to the Picard functor. 
\end{nota}
\begin{nota}{\it References.}
In the literature, there exist several  constructions of compactified Picard schemes
(see for example  \cite{I}, \cite{dsouza}, \cite{OS}, \cite{AK}, \cite{simpson} 
\cite{caporaso}, \cite{pandha}, \cite{ed01})
which differ from one another in various  aspects, such as the  functorial interpretation.
A survey may be found in \cite{alex}.

Our  goal here is to classify   them according to their relation with the N\'eron model; see  \ref{NDtype}.
In doing so  we shall not go through the details concerning various constructions, but we shall 
focus on    their
 formal properties. 
\end{nota}
We first describe   an apparently  simple, yet  challenging case.

\begin{example} {\it }
\label{prov} The following is a special case
of Example~\ref{vine}, whose notation we continue to use. Let $X=\vine=C_1\cup C_2$  with $\dd =1$. So $X$ is a
curve of compact type whose generalized jacobian is projective (cf. \ref{genjac}):
$$
J(X)=\Pic^{(0,0)}X\cong\Pic^0C_1\times \Pic^0C_2=J(C_1)\times J(C_2).
$$
If $f:\X\to B$ is a regular smoothing of $X$, the Picard scheme $\Pic_f^d\to B$ is   not separated
(cf. \ref{relpic}).
Let us look for 
a separated,  even projective, model for $\Pic^d\gen$; in other words, let us look for a compactified Picard
scheme.  This does not seem too hard: every connected component
of $\Pic^dX$ is projective.
Recall also (see example~\ref{dcgvine}) that   $\dcg$ consists of only one element, so we  expect  a separated model
of $\Pic^d\gen$ to have  only one irreducible component.

If $d=0$ there is no problem: it suffices to take the generalized jacobian 
$J_f\to B$ which is certainly projective. In doing so we have made a choice of multidegree for the closed fiber,
namely the multidegree  $(0,0)$. 

If $d\neq 0$   there does not seem to be a  canonical choice of multidegree.
Indeed,   consider the case $d=g-1$ (which will turn out to be quite   interesting in the sequel).
Since $g=g_1+g_2$,   there exists no natural choice of a multidegree of total degree $g-1$.
The lesson we should draw from this is that
a compactified degree-$d$ Picard scheme for $X$ may fail to contain a subset corresponding to line bundles 
of degree $d$ on
$X$ itself  (this does not happen if $X$ is irreducible; see \ref{N=Dirr}).

In the present example, we shall see that there exists a canonical compactified Picard scheme $\PXgb$
 endowed with a canonical isomorphism with $\Pic^{(g_1-1,g_2-1,1)}\hX$, where $\hX=C_1\cup C_2\cup C_3$
is the  curve  obtained by ``blowing-up" $X$ at the node and calling the  ``exceptional" component 
$C_3$. 

An equivalent description of $\PXgb$
is as the moduli space of   Euler characteristic 0, rank-1, torsion-free sheaves on $X$ that are not locally free at
the node.

In either case, $\PXgb$ does not parametrize any set of line bundles of degree $g-1$ on $X$.
The  above description     yields an isomorphism 
$$
\pXgb \cong \Pic^{(g_1-1,g_2-1)}X,
$$
so, an interpretation in terms of line bundles  of degree $g-2$
rather than  $g-1$ (see \ref{descg-1} for more details).
\end{example}

\begin{nota}{\it N\'eron models and compactified Picard schemes.}
\label{types}
Let  $X$ be any curve and $f:\X \to B$ a regular smoothing for $X$.
Let us now consider a fixed  compactifed Picard scheme  $\pi_d:\pfb\to B$ 
and its closed fiber $\PXb$   (cf. \ref{comppic}). 
Recall that we call $\pf\to B$ the smooth locus of
$\pi_d$. We  apply the N\'eron mapping property to obtain  
  a canonical $B$-morphism 
\begin{equation}
\label{nf}
n_f: P^d_f\to \nf
\end{equation}
extending the indentity map from the generic fiber of $\pi_d$ to the generic fiber of  $\nf\to B$.
\end{nota}
\begin{defi}
\label{NDtype} 
We say that $\PXb$, or $\pfb$, is   of {\em N-type} (or of {\em N\'eron type}) if  the map $n_f:\pf \to \nf$ 
 is an isomorphism.
We say that $\PXb$, or $\pfb$, is of {\em D-type} (or of {\em Degeneration type})  otherwise.
\end{defi}

\begin{nota}{\it Irreducible curves.}
\label{N=Dirr}
If $X$ is   irreducible,  up to isomorphisms there exists in the literature 
a unique compactified Picard
scheme
$\PXb$, which is of
N-type and does not depend on $d$. 
$\PXb$ is irreducible and reduced; it is singular unless $X$ is smooth.
It is also true that $\Pic^dX$ is naturally isomorphic to an open subset of $\PXb$
(the first such constructions are \cite{I}, \cite{MM}, \cite{dsouza} and \cite {AK}).
\end{nota}

As we said, if $X$ is reducible and not tree-like, then the structure and the type of $\PXb$ varies. 
\begin{example}
Let us  consider a curve $\vine=C_1\cup C_2$ as in Example~\ref{vine}, and assume $\dd \geq 2$.

The   fact that there exist different, non isomorphic, compactified jacobians for $\vine$ was first
discovered by Oda and Seshadri in their fundamental paper \cite{OS}, 
the first paper   providing a construction of a compactified jacobian for any reducible nodal fixed curve.
They proved that, depending on a certain choice of polarization,
there exists s a compactified jacobian of  N-type, hence having $\dd$ irreducible components,
or a compactified jacobian of  D-type, having $\dd -1$ irreducible components (see \cite{OS} Chapter II, section
13).

The same was later found in \cite{caporaso} in a different framework. 
There is no polarization involved in this construction, and the structure of $\PXb$ depends solely
on the degree
$d$. For every curve $\vine$ both types of compactification  appear, and they are isomorphic 
to the ones of
\cite{OS}. More precisely, if $d=g-1$ then for every $\dd\geq 2$,
we have that $\PXbd$ has $\dd-1$ components
(hence it is of D-type). If $d=0$ then $\PXbd$ is of N-type  (resp. 
of D-type with $\dd -1$ components) if $\dd$ is odd (resp if $\dd $ is even). 
See Theorem~\ref{factMR} for values of $d$ for which $\PXbd$ is of N-type.
\end{example}

\begin{nota}{\it Reducible  curves.}
\label{redND}
Once we
choose the specific construction (so that we choose a  $\PXb$ for every $X$ and $d$) 
 the following problem arises naturally: 
for a fixed $d$ classify all  curves for which $\PXb$ is of N\'eron type.
A satisfactory answer to this question is known if $d=g-1$ (and the answer is: only the obvious ones; see
Theorem~\ref{factg-1} below), and in    a few other cases. 
A construction that is rather well understood  is that  of \cite{caporaso}
or
\cite{pandha} (the two  are naturally isomorphic);
here is what is known.
\begin{thm}
\label{factMR}
For every $d$ there exists an open substack $\mgbst  ^d$ of $\mgbst$, whose  moduli scheme $\mgbar  ^d$
we call
 the locus of {\em $d$-general}
stable curves, such that the following holds.

\begin{enumerate}[(i)]
\item
There exists a Deligne-Mumford stack $\pdbstn$ with a strongly representable morphism
onto $\mgbst  ^d$, such that for every $X\in \mgbar  ^d$ and every regular smoothing $f:\X\to B$ of $X$,
the  base change $\pdbstn\times _{\mgbst}B\to B$ is a compactified Picard scheme   of N\'eron-type.
\item
 $\mgbar  ^d$ is open in $\mgbar$ and contains  the locus of tree-like curves.
\item
\label{MR}
$\mgbst  ^d=\mgbst $ if and only if $(d-g+1,2g-2)=1$.
\item
\label{g-1gen}
If $d=g-1$ then $\mgbar  ^d$ equals  the locus of tree-like curves.
\end{enumerate}
\end{thm}
See  \cite{cner} for the case described in part (\ref{MR}) 
and \cite{melo} for the remaining cases;  in the latter paper there is also a
complete  description of $\mgbar  ^d$, based on the combinatorics of the curves.

By \cite{alex} section 1, Theorem~\ref{factMR}, which is proved using the \cite{caporaso} construction, applies also
for the compactifications of \cite{OS} and \cite{simpson} with respect to the canonical polarization.

Finally,  the construction  of \cite{ed01}  (which
concerns    curves with singularities more general than nodes)
 has been studied  in \cite{busoneroest}. It is
there shown that the compactified Picard schemes
called $J^{\sigma}_{\mathcal E}$  are of N-type for all regular smoothings
(endowed with a section $\sigma$ and a polarization 
${\mathcal E}$, needed for the construction).

\end{nota}
\subsection{The case $d=g-1$}

\

\

As we mentioned, 
 the case $d=g-1$   has been studied  closely by various authors and it is thus much better understood.
The following is a   summary of known results:
\begin{thm}
\label{factg-1}
For every  curve $X$ of genus $g\geq 2$ the following facts hold.
\begin{enumerate}[(i)]
\item
\label{all=}
The
compactified degree-$(g-1)$ Picard schemes   constructed in
\cite{OS},
\cite{simpson} and    \cite{caporaso} are all isomorphic 
 to a projective, reduced scheme,  denoted $\PXgb$ from now on.
\item\label{theta}
$\PXgb$ possesses a theta divisor $\Theta(X)$
which is Cartier and
ample. If  $X$ is smooth, $\Theta(X)$ coincides with the classical theta divisor
i.e.  $\Theta(X)=W _{g-1}(X)$ in the standard notation (see \cite{ACGH}).
\item
\label{ssp} The pair $(\PXgb, \Theta(X))$ is a semiabelic stable pair in the sense of \cite{alex1}.
\item
\label{D}
If $X$ is  not   tree-like,    $\PXgb$ has less than $c(X)$ components (hence it is
of
$D$-type).
\end{enumerate} 
\end{thm}
Parts (\ref{all=}) and (\ref{ssp}) are due to V. Alexeev,  \cite{alex}.
Part (\ref{theta}) is due to A. Soucaris \cite{soucaris} and E. Esteves \cite{esttheta} in case $X$ is irreducible,
and  to Alexeev if is $X$ reducible in which case the work of A. Beauville in
\cite{beau} plays a key role;
 see  \cite{alex} for  details.
For part (\ref{D})   see  \cite{cner} Section 4.

\begin{nota}{\it The canonical compactified jacobian in degree $g-1$.}
\label{g-1can}
As we said, the above results, especially part (\ref{ssp}),  lead us to regard such compactification $\PXgb$
as canonical. 
Throughout   this paper, the compactified Picard schemes 
$\PXgb$  and $\ov{P_f^{g-1}}$ will be  the ones of Theorem~\ref{factg-1}.
We shall give it a more explicit
description in the sequel.

Having Theorem~\ref{factMR} in mind we now ask:
how do these  Picard schemes $\PXgb$ glue together over $\mgbst$,
as $X$ varies among all stable curves of genus $g$?
 Of course
Theorem~\ref{factMR} does not tell us much.
The following result  answering this question is due to   M. Melo; see  \cite{melo}.
\end{nota}
\begin{prop}
\label{MM}
There exists an Artin stack $\overline{{\mathcal{P}}_{g-1,g}}$ with a non representable morphism to 
$\mgbst$, such that for
every $X\in \mgbar$ and every smoothing $f:\X\to B$ of it, 
the stack $\overline{{\mathcal{P}}_{g-1,g}}\times_{\mgbst}B$
admits a canonical  proper $B$-morphism onto $\pfgb$.
\end{prop}

\begin{remark}
\label{g-1}
In particular, the  stack $\overline{{\mathcal{P}}_{g-1,g}}$ is not  Deligne-Mumford. 
We like to interpret this
phenomenon as a reflection of the fact that $\PXgb$, being of D-type, does not have the best moduli properties that
one may hope for. Recall in fact that
$\PXgb$  has   fewer components than
the N\'eron model (by \ref{factg-1} (\ref{D})). This  tells us that some multidegree classes are not ``finely"
represented by points in $\PXgb$ (see Example~\ref{vine1}).

However, if we restrict our attention to the moduli scheme of automorphism free stable curves, $\mgbar^0$
(so that there is a universal family   $\ov{{\mathcal C}_g}\to \mgbar^0$),
then there does exist a scheme 
\begin{equation}
\label{g-1sch}
\overline{P_{g-1,g}}\la \mgbar^0
\end{equation}
whose fiber over every curve $X$ is the (canonical) $\PXgb$.
\end{remark}

\begin{example}
\label{vine1}
Consider the  curve 
$
\vine=C_1\cup C_2
$
and assume that $\dd\geq 2$
(notation in Example~\ref{vine}).
Then $\PXgb$
 has $\dd -1$ irreducible components each of which contains a copy of
$J(X)$ as a dense open subset.
So, we seem to have lost a multidegree class (cf. example~\ref{dcgvine})!

What actually happens
 is that there is one   multidegree class,
call it $\mu_0\in \dcg^{g-1}$, such that line bundles having multidegree  of class $\mu_0$ are represented
by points in the boundary of $\PXgb$. 
Furthermore, different such line bundles get identified.

The simplest case when that happens is  $\dd =2$. 
$\PXgb$ is thus irreducible,  and it turns out to contain a dense open subset (equal to its smooth locus)
naturally identified to
$\Pic^{(g_1,g_2)}X$
(see Proposition~\ref{Pstr} and Definition~\ref{bal}). We shall henceforth identify
the smooth locus of $\PXgb$ with $\Pic^{(g_1,g_2)}X$.

The class $\mu_0$ 
(defined above) is  thus $\mu_0=[(g_1-1,g_2+1)]=[(g_1-1-2n,g_2+1+2n)]$, $n\in \Z$.
Now the boundary of $\PXgb$ is an irreducible $(g-1)$-dimensional closed subscheme
which is isomorphic to the jacobian of the normalization of $X$. More precisely,
we have that the boundary has a canonical isomorphism

\begin{equation}
\label{boundary}
\PXgb\smallsetminus \Pic^{(g_1,g_2)}X \cong \Pic^{(g_1-1,g_2-1)}\Xn=\Pic^{g_1-1}C_1\times \Pic^{g_2-1}C_2
\end{equation}
where $\Xn=C_1\coprod C_2$ is the normalization of $X$.

Let $L\in \Pic^{(g_1-1,g_2+1)}X$, pick a regular smoothing $f:\X \to B$ of $X$ such that there exists 
a line bundle $\L$ on $\X$  restricting to $L$ on the closed fiber.
Then there is a  map $\phi:B\to \ov{P^{g-1}_f}$ such that $\phi(\Spec K)=[\L_{\gen}]$
(of course, $\phi$ is regular as $B$ is a smooth curve and $\ov{P^{g-1}_f}$ is projective).
By what we claimed before, $\phi$ must map the closed point of $B$ to a boundary point of $\pXgb$;
which point?

The answer is, using (\ref{boundary}), the point of the line bundle
$(L_1,L_2(-p_2-q_2))$ on
$\Xn$, where $L_i$ denotes the restriction of $L$ to $C_i$, and $p_2,q_2$ are the branches of the nodes of $X$ lying
in
$C_2$ (see also \ref{funct}).

We conclude by noticing that this gives a map from $\Pic^{(g_1-1,g_2+1)}X$ to the boundary
of $\PXgb$, which is surjective and has one-dimensional fibers.
More details will be in \ref{funct}.
\end{example}

\begin{nota}{\it A stratification of   $\PXgb$.}
\label{descg-1} 
Recall  a well known
\begin{defi}
\label{bal} Let $\md \in \Z^{\gamma}$ be such that $|\md|=g-1$. Then $\md$ is {\it semistable} (resp. 
{\it stable})
on
$X$ if for every connected    subcurve $Z\subsetneq X$ we have $d_Z\geq p_a(Z)-1$ (resp. $d_Z> p_a(Z)-1$),
where $d_Z=\sum_{C_i\subset Z}d_i$ denotes the total degree of the restriction of $\md$ to $Z$
and $p_a(Z)$   the arithmetic genus of $Z$.

If $Y$ is a nodal, disconnected curve, we say that a multidegree is semistable or stable 
if it is so on every connected component of $Y$.

We denote by $\BXs$ (resp. by $\BX$) the set of semistable  (resp. stable) multidegrees on $X$.
\end{defi}

If $X$ is irreducible, then of course $\BXs=\BX=\{g-1\}$.

It is easy to check that $\BXs$ is finite and not empty, whereas $\BX$ may be empty
(see example~\ref{Pvine}).

It  is well known that every class in $\dcg^{g-1}$ has some semistable representative.

\begin{example}
\label{Pvine}
If $\vine$ is the vine curve of example~\ref{vine}, then
$$
\#\Sigma^{ss}(\vine) = \dd+1 \  \text{ and }\   \  \#\Sigma(\vine) = \dd-1.
$$ More precisely
$$
\Sigma^{ss}(\vine) =\{(g_1-1,g_2-1+\dd),\ldots,(g_1-1+\dd,g_2-1)\}
$$
and $\Sigma(\vine) =\{(g_1,g_2-2+\dd),\ldots,(g_1-2+\dd,g_2)\}$. In particular,
$\Sigma(\vine) $  {\it  is empty if and only if }
 $\dd=1$.

If $\Xn=C_1\coprod C_2$ is the normalization of $X$ (so that $\Xn$ is disconnected),
one easily checks that 
\begin{equation}
\label{desc}
\Sigma^{ss}(C_1\coprod C_2)=\Sigma(C_1\coprod C_2)=\{(g_1-1,g_2-1)\}.
\end{equation}
\end{example}
Here  is a  coincise way of describing  $\PXgb$
(for the proof, one may choose any of the constructions \cite{OS}, \cite{caporaso}, \cite{simpson}
and use their equivalence  established in \cite{alex}  (cf. Theorem~\ref{factg-1}).
\end{nota}

\begin{prop}
\label{Pstr}
The points of $\PXgb$ bijectively parametrize all line bundles having stable multidegree 
on all partial normalizations of $X$ (including $X$ itself).
\end{prop} 

More precisely, for every subset $S\subset \sing$ consider the partial normalization 
$\XS\to X$ of $X$ at exactly
$S$. Consider now $\Sigma (\XS)$; if it is nonempty let $\md $ be a multidegree in it (note  that $|\md|=p_a(\XS)-1$)
 and consider $\Pic^{\md}\XS$. Then there exists a canonical injective morphism
$\eSd:\Pic^{\md}\XS\ha \PXgb$ whose image we denote
$$
\PSd:=\eSd(\Pic^{\md}\XS)\subset \PXgb.
$$
The sets $\PSd$ form a stratification of $\PXgb$ in disjoint strata, i.e. we have
$
\pXgb =\coprod_{S \subseteq\sing \  \md \in \BXS} \PSd.
$
More details can be found in \cite{Ctheta}.

The explicit description of $\eSd$ depends on how we choose to describe 
$\PXgb$ functorially
(whether we use rank-1 torsion free sheaves, line bundles on semistable curves, cell
decompositions...). In the next example we will use line bundles on semistable curves.

\begin{example}
\label{strvine} 
Let $X=\vine=C_1\cup C_2$. 
If we interpret $\PXgb$ as parametrizing line bundles on blow-ups
of $X$, then the map $\eSd$ is described as follows.
If $S=\emptyset$ and $\md$ is a stable multidegree on $S$ (in particular, $\dd \geq 2$, by example~\ref{Pvine}),
then $\eSd$ is the identity map.

If $S=\{n\}$ is one node, let $\XS$ be the normalization of $X$ at $n$ and $\hXS=\XS\cup C_3$
be the curve obtained by joining the two branches over $n$ by a smooth rational curve $C_3$.
Let $\md \in \BXS$, so that $|\md|=p_a(\XS)-1=g-2$.
Now $\eSd$ maps a line bundle $L\in \Pic^{\md}\XS$ to the (unique) point in  $\PSd\subset \PXgb$
corresponding to line bundles on $\hXS$ whose restriction to $\XS$ is $L$ and whose restriction to
$C_3\cong \pr{1}$ is $\O(1)$.

Finally, if $X=\sing$ then $\XS=\Xn$ is the normalization of $X$. Using Example~\ref{Pvine}
we see that $\eSd=\epsilon^{(g_1-1,g_2-1)}_{\sing}$. With a procedure analogous to the previous case,
 we get that the
smallest stratum of $\PXgb$ is isomorphic to $\Pic^{g_1-1}C_1\times \Pic^{g_2-1}C_2$.
Call  $\hX $ the (connected, nodal) curve obtained by blowing up every node of $X$
so that $\hX$ is the union of $\Xn$ with $\delta $ copies of $\pr{1}$, one for each node.
Now every point $\ell$ of this stratum corresponds to the set of line bundles
on $\hX$  whose restriction to $\Xn$ is a fixed line bundle $L$
of multidegree $(g_1-1,g_2-1)$, and 
whose restriction to each of the remaining components is $\O(1)$. Thus $\ell = \eSd(L)$.

\end{example}
\begin{nota}{\it Functorial properties of $\pfgb$.}
\label{funct}
Fix the curve $X$ and a regular smoothing $f$ for it.
We are going to illustrate some moduli properties of
$\pfgb$,  using the same notation as for the moduli property of $\Pic_f$ (cf. \ref{relpic}).
For every $B$-scheme $T\to B$  and  every line bundle $\L$ on $\X_T$,
such that for every $t\in T$ the restriction $\L_{f_T^{-1}(t)}$ has semistable multidegree, there exists a canonical
morphism
$$
\ov{\mu} _{\L}:T\la \pfgb
$$
having the following properties.
The restriction of $\ov{\mu}_{\L}$ over the generic point is the usual
moduli map to the Picard variety $\Pic^{g-1}\gen$. If $t$ lies over the closed point of $B$, then obviously
$f_T^{-1}(t)=X$ and we are assuming that $\mdeg \L_{|f_T^{-1}(t)}\in \BXs$.
We denote
$\ov{\mu}_{\L}(t)=[\L_{|f_T^{-1}(t)}]\in \PXgb$. What
is this class, in terms on the description given in  \ref{descg-1}?

Call $L=\L_{|f_T^{-1}(t)}$ and $\md$ its multidegree. If $\md$ is stable, then there is no ambiguity:
by Proposition~\ref{Pstr} there is a  point $[L]$ in $\PXgb$ corresponding to $L$.

If $\md$ is strictly semistable some cumbersome notation  is needed for an arbitrary curve.
Therefore, to give a more efficient explanation,
we shall  precisely describe only the case of a curve with two components.

So, let $X=\vine=C_1\cup C_2$, recall (cf. example~\ref{Pvine}) that there are exactly two
strictly semistable multidegrees, and they are  equivalent;
namely 
$$
(g_1-1, g_2+\dd -1)\equiv (g_1+\dd-1, g_2 -1).
$$
Call $L_1$ and $L_2$ the restrictions of $L$ to $C_1$ and $C_2$.
Denote $\{ p_1,\ldots,p_{\dd}\}\subset C_1$ (resp. $\{ q_1,\ldots,q_{\dd}\}\subset C_2$)
the $\dd$ points of $C_1$ (resp. of $C_2$) lying over the nodes of $X$.

If $\mdeg L=(g_1-1, g_2+\dd -1)$ then $\ov{\mu}_{\L}(t)$ is  the point in the stratum
$P_{\sing}^{(g_1-1, g_2-1)}\cong \Pic^{(g_1-1, g_2-1)}\Xn$ given by
$$
[(L_1, L_2(-\sum_{i=1}^\dd q_i)]
$$
(see \ref{desc} and example \ref{Pvine}).
If instead $\mdeg L=(g_1+\dd -1, g_2-1)$,
then
$$
\ov{\mu}_{\L}(t)=[(L_1(-\sum_{i=1}^\dd p_i), L_2].
$$
\end{nota}
\section{Abel maps and Theta divisors}
\subsection{Preliminary analysis}

\begin{nota}{\it The smooth case.}
\label{smooth}
Let $C$ be a    smooth curve. For every integer $d\geq 1$ the $d$-th Abel map
is defined as follows
$$
\alpha^d_C:C^d\la \Pic^dC ;\  \  (p_1,\ldots,p_d)\mapsto \O_C(\sum p_i);
$$
$\alpha^d_C$  is a regular map.
It is well known that  Abel maps are defined more generally for any family of smooth curves over any scheme.
Moreover the image of the $d$-th Abel map in $\Pic^dC$ is equal to the variety  
$W_d(C):=\{L\in \Pic^dC: h^0(C,L)\neq 0\}$, and that
\begin{equation}
\label{AW}
\im \alpha^d_C = W_d(C)  \  \  \text{ and }    \dim W_d(C) =\min \{d,g\}.
\end{equation}

The situation is particularly interesting if $1\leq d\leq g$; then $W_d(C)$ is a proper subvariety of $\Pic^dC$.
Moreover, for any nonnegative integer $r$,
 the loci in $W_d(C)$ where the fiber dimension of the $d$-th Abel map is
at least   $r$
are the Brill-Noether varieties $W^r_d(C)$ (so that $W_d(C)=W^0_d(C)$).
The geometry  of line bundles and linear series on a smooth curve $C$ is encoded in the varieties $W^r_d(C)$;
see \cite{ACGH} for the general theory, from the time of Riemann to
the twentieth  century.

How does this picture extend  to singular curves?
\end{nota}

\begin{nota}{\it Naive approach for singular curves.}
Fix a degree $d\geq 1$ and a  curve $X$. One may define a rational map
\begin{equation}
\label{naive}
\begin{array}{lccr}
\widetilde{\alpha}^d_X : & X^d&\da &\Pic^dX \\
&(p_1,\ldots, p_d) &\mapsto &\O_X(\sum_1^dp_i)
\end{array}
\end{equation}
which is regular if all the $p_i$ are nonsingular points of $X$. The above definition is the simple 
minded
extention
of the  smooth curve case, and it turns out   to be non-satisfactory,
unless $X$ is irreducible (see Proposition~\ref{path}). To be more precise, denote by $X=\cup_{i=1}^{\gamma}C_i$ the
irreducible component decomposition of $X$ and  for any $\md=(d_1,\ldots,d_{\gamma})\in \Z^{\gamma}$ such that $|\md|=d$
and $\md \geq 0$ (i.e. $d_i\geq 0,\  \forall i$),  set
\begin{equation}
\label{Xmd}
\Xmd :=C_1^{d_1}\times \ldots\times C_{\gamma}^{d_{\gamma}};
\end{equation}
more generally, for any permutation $\sigma$ of the set $\{1,\ldots, \gamma\}$, let 
\begin{equation}
\label{Xmds}\Xmds :=C_{\sigma(1)}^{d_1}\times \ldots\times
C_{\sigma(\g)}^{d_{\gamma}}.
\end{equation}
Thus the  $\Xmds$ are the irreducible components of $X^d$.
If $\sigma$ is the identity we often omit it (as in (\ref{Xmd})). Now let
\begin{equation}
\label{amd}
\amd:\Xmd \  \da \picX{\md}
\end{equation}
(respectively, $\amds:\Xmds \da \picX{\md}$)
be the restriction of  $\widetilde{\alpha}^d_X$ to $\Xmd$
(respectively, to $\Xmds$). These maps are of course defined only if every
$d_i$ is nonnegative. To simplify matters, whenever they are   not defined, we shall set
$\im \amds =\emptyset$.

We   define for any curve $X$ and any multidegree $\md$
\begin{equation}
\label{Wmd}
\Wmd:=\{L\in \Pic^{\md}X: h^0(X,L)\geq 1\}\subset \picX{\md}
\end{equation}
and $W_d(X)=\coprod_{|\md|=d}\Wmd$. In analogy with \ref{smooth}
 we   ask what is the relation between $\amd$ and $\Wmd$.
Here is where  the first type of pathologies appears.
The following statement is interesting only for a reducible curve $X$;
we leave it to the reader to find the analogue for $\amds$.
For any component $C_i$ of $X$ denote by $g_i$ its arithmetic genus and by $\dd_i=\#(C_i\cap \ov{X\smallsetminus
C_i})$.
\end{nota}
\begin{prop}
\label{path}
Let $X$ be a curve of genus $g\geq 2$  and $\md \in \Z^{\gamma}$ with $|\md|=d\geq 1$.
\begin{enumerate}[(1)]
\item
\label{pathd} If there exists an   $i$ such that $d_i\geq g_i +\dd_i $, then
$\dim \Wmd =g $ and    $\dim \im\amd\leq
d-1. $
\item
\label{pathg-1}
Assume that   $|\md|=g-1$.  
\end{enumerate}
\begin{enumerate}[(a)]
\item
\label{nss}
If $\md$ is not semistable, then $\dim \Wmd =g$ and $\dim \im\amd\leq g-2$.
\item
\label{ss}
If $\md$ is semistable, then $\dim \im\amd =\dim \Wmd=g-1$.
\end{enumerate}
\end{prop} 
\begin{proof}
Part (\ref{nss}) is a special case of 
Part (\ref{pathd}).
Part (\ref{ss})
 combines some results of Beauville,  namely Lemma 2.1 and Prop. 2.2 in \cite{beau}, with some of
\cite{Ctheta} (Prop. 3.6).

Let us prove part (\ref{pathd}).
Assume that $d_i\geq g_i+\dd_i$; then for every $L_i\in \Pic^{d_i}C_i$ we have
$h^0(C_i,L_i)\geq \dd_i+1$.
Therefore for every $L\in \picX{\md}$, $L$ admits some global section that does not
vanish on $C_i$.
Hence $h^0(X,L)\geq 1$ and $\Wmd=\picX{\md}$. This proves the first statement of part
(\ref{pathd}). 

For the second, suppose that $\amd$ is defined (i.e. that $\md \geq 0$).
Call $C_i'=\ov{X\smallsetminus C_i}$ the complementary curve of $C_i$
and let
$\nu_i:C_i\coprod C_i'\la X$ be the normalization of $X$ at $C_i\cap C_i'$.
Let $L\in \im\amd$ and set  $M:=\nu_i^*L$ , so that $M$ is determined by a pair of line bundles
$L_i\in \Pic C_i$ and $L_i'\in \Pic C_i'$.
We have that $h^0(C_i,L_i)\geq \dd_i+1$ (by what we said before) and $h^0(C_i',L_i')\geq 1$
because $L\in \im \amd$. Therefore $h^0(M)\geq \dd_i+1+1=\dd_i+2$; this  implies that
$$h^0(X,L)\geq h^0(M)-\dd_i\geq 2.$$

We conclude that the fibers of the map $\amd$ have dimension at least equal to $1$, and hence that
$\dim \im \amd \leq \dim \Xmd-1 = d-1$ as wanted.
\end{proof}

\subsection{Abel maps in degree $g-1$}

\begin{nota}{\it Defining Abel maps by specialization.}
 Proposition~\ref{path}  indicates that dealing with Abel maps of reducible curves is  a delicate
matter and more sophisticated tools may be needed.
Therefore,  we shall  now  define Abel maps for a curve $X$  after a smoothing $f$ of $X$ is given,
and using compactified Picard
schemes. In other words, 
  we shall add  some variational data,  defining Abel maps for a so called   ``enriched curve", i.e. a pair
 $(X,f)$. 
Furthermore we shall use a compactified Picard scheme as the image space of our Abel maps.
This will give us a better behaved object, yet one that may (and will) depend on the the choice of the smoothing,
and on the choice of the compactified Picard scheme.
\end{nota}
\begin{nota}{\it Abel maps in degree $g-1$.}
\label{abelg-1}
We first deal with the case $d=g-1$, so that
 we have a canonical choice for the compactified Picard scheme
(see \ref{g-1can}). 

Let $X$ and $f:\X\to B$ be a  curve and a  regular smoothing for it. Denote by
$\X^{d}=\X\times_B\ldots\times_B \X$ the $d$-th fibered power over $B$.
Classically, the $d$-th Abel map is the moduli map associated to the universal effective
divisor on $\X^{d}\times _B\X$ (see below). We shall approach the problem in the same way.
So consider the natural projection
$$
\X^g=\X^{g-1}\times _B\X \stackrel{\pi}{\la }\X^{g-1}
$$
(above and throughout this section,
 all schemes, maps and products are over $B$). The map $\pi$ has $g-1$ tautological rational sections
$\sigma_i(p_1,\ldots,p_{g-1})=(p_1,\ldots,p_{g-1};p_i)$ which determine a line bundle $\mathcal E$
on the smooth locus of
$\X^g$, i.e. introducing the  (Weil) divisor
$$
E=\sum_1^{g-1}\ov{\sigma_i(\X^{g-1})}
$$
we define
$$
\E=\O_{\X^g}(E)
$$
which is  locally free on an   open subset of $\X^g$.
Now, the $\pi$-relative   multidegree of $\E$ varies with the irreducible components of $X^{g-1}$.
Indeed, with the notation of
(\ref{Xmd}), let $\Xmd\subset X^{g-1}\subset \X^{g-1}$ be an irreducible component (so that $\md \geq 0$),
then   for a generic point $t\in \Xmd$ we have
$$
\mdeg \E_{|\pi^{-1}(t)}=\md.
$$
More generally, if $t\in \Xmds$  we have 
$\mdeg \E_{|\pi^{-1}(t)}=(d_{\sigma^{-1}(1)},\ldots, d_{\sigma^{-1}(\gamma)})$.
 Observe  that the restriction
of $\E$ over $\Spec K$
is the so-called universal effective divisor on $\gen^{g-1}\times \gen$, whose moduli map $\gen^{g-1}\to
\Pic^{g-1}\gen$ is the $g-1$-th Abel map of $\gen$. So, 
we would like to  complete this, associating to $\E$ a map
$\X^{g-1}\da \pfgb$.
With the functorial description of $\pfgb$ in mind (see
\ref{funct}), the question we need an answer for  is: is $\md$ semistable?

To better explain how to proceed we   concentrate on our leading example.
\end{nota}
\begin{example}
\label{abelvine}
Assume that $X=\vine$,
so that we have a simple description of $\BXs$ and hence of $\PXgb$.
The irreducible component decomposition of $X^{g-1}$ is
$$
X^{g-1}=\bigcup_{l=1}^{g-1}\bigr( C_1^l\times C_2^{g-1-l}\cup C_2^l\times C_1^{g-1-l}\bigr)
$$
so that if $t\in C_1^l\times C_2^{g-1-l}$, then $\deg \E_{|\pi^{-1}(t)}=(l, g-1-l)$
while if $t\in C_2^l\times C_1^{g-1-l}$, then $\deg \E_{|\pi^{-1}(t)}=(g-1-l,l)$.
Now we claim that for every $l=0,\ldots,g-1$ there exists an integer $a(l)\in \Z$ such that
\begin{equation}
\label{adj} 
(l, g-1-l)+(-a(l)\dd, a(l)\dd)\in \BXs.
\end{equation}
Indeed, $(-a(l)\dd, a(l)\dd)\in \Lambda_X$, so $(l, g-1-l)+(-a(l)\dd, a(l)\dd)\equiv (l, g-1-l)$;
as every multidegree class has a semistable representative 
the claim is proved.

If $(l, g-1-l)$
is semistable we  choose $a(l)=0$.
Note that $a(l)$ may not be unique, but this will turn out to be irrelevant. 
Indeed, $a(l)$ is not unique if and only if $(l, g-1-l)+(-a\dd,+a\dd)$ is strictly semistable.
From the description given in \ref{funct} one sees that the choice of $a(l)$ plays no role.

Now we define
\begin{equation}
\label{Lg}
\Lg=\E\otimes \O_{\X^g}\Bigr(\sum_{l=1}^{g-1}a(l)\bigr( C_1^l\times C_2^{g-1-l}\times C_1
+C_2^{g-1-l}\times C_1^l\times C_1\bigl)\Bigl)
\end{equation}
which is locally free over the smooth locus of $\X^{g-1}$.
Let $t$ be a smooth point of $X^{g-1}$;
by construction, if either  $t\in C_1^l\times C_2^{g-1-l}$
or $t\in  C_2^{g-1-l}\times  C_1^l$  we have that
$$
\mdeg (\Lg)_{|\pi^{-1}(t)}=(l-a(l)\dd,g-1-l+a(l)\dd)
$$
which is semistable, by (\ref{adj}).
Therefore there exists a canonical rational map 
$
\ov{\mu}_{\Lg}:\X^{g-1}\da \pfgb
$  (see \ref{funct}). 

We denote $\ag:=\ov{\mu}_{\Lg}$ and call it the $g-1$-th Abel map associated to $f$.
The restriction of $\ag$ to the closed fiber is the rational map
\begin{equation}
\label{agX}
\agX:X^{g-1}\da \PXgb.
\end{equation}
This is the definition we were aiming at; so we call $\agX$ the  $g-1$-th Abel map of $X$ associated to $f$.
By construction, 
$\agX$ is regular at $(p_1,\ldots, p_{g-1})$
for every $p_1,\ldots, p_{g-1}\in X\smallsetminus \sing$.
\end{example}
We just gave the definition in the special case of a curve with two components.
The general case can be dealt with using the very same procedure,  paying quite a   price in terms of notation.
Rather than going through this, we prefer to deal with a problem that arises immediately.

\begin{nota}{\it Naturality.}
\label{g-1nat} We now consider the following question:

{\it Does  $\agX:X^{g-1}\da \PXgb$ depend on the choice of the smoothing $f$?}

We shall say that $\agX$ is {\it natural} if it is independent on the choice of $f$, i.e. 
if for every regular smoothings $f$ and $f'$ of $X$, we have
$$
\agX=\alpha^{g-1}_{f',X}.
$$

We have defined $\agX$ only for a curve $X=\vine$, so we shall focus on this case,
which is already interesting. See \cite{busonero} for the general result, valid for all stable curves.
We use the notation of Example~\ref{vine}, we have
\begin{prop}
\label{natg-1}
Let $X=\vine$.
If $\dd =1$ then $\agX$ is natural.

Assume $\dd \geq 2$, then
$\agX$ is natural if and only if 
$$
\dd \geq g-1 \text{ and }\  \{ g_1,g_2\}\neq \{0,2\}.
$$
Equivalently, $\agX$ is natural if and only if  $g_i\leq 1$ for $i=1,2$.
\end{prop}
\begin{proof} The map $\ag$ is defined as the moduli map $\ov{\mu}_{\Lg}$, and 
$\Lg$ is a so called ``twist" of $\E=\O(E)$
(i.e. $\Lg$ and $\E$ differ only over the closed point of $B$).
It is clear that the restriction of $\E$ to the fibers over $X^{g-1}$ is independent of $f$
(indeed $\E_{|\pi^{-1}(p_1,\ldots, p_{g-1})}= \O_X(p_1+\ldots, p_{g-1})$).

Now consider the other factor ${\mathcal T}:=\O_{\X^g}\Bigr(\sum_{l=0}^{g-1}a(l)\bigr( C_1^l\times C_2^{g-1-l}\times C_1
+C_2^{g-1-l}\times C_1^l\times C_1\bigl)\Bigl)=\Lg\otimes  {\mathcal E}^{-1}$
of (\ref{Lg}). The restriction of  $\mathcal T$ to the fibers of $\pi$ is a twister.
Now recall that if $\dd\geq 2$
a nontrivial twister on $\vine$ depends on the choice of $\X$ (see remark~\ref{dep}).
On the other hand if $\dd =1$, a twister on $X$ is uniquely determined by its multidegree.
Hence if $\dd =1$ the map $\agX$ does not depend on $f$.

Assume from now on $\dd \geq 2$. By what we said, $\agX$ is natural iff the map is not twisted iff 
${\mathcal T}=0$.
Now, ${\mathcal T}=0$ iff the multidegree $(l,g-1-l)$ is semistable for every $l=0,\ldots,g-1$
(by \ref{abelg-1}).

Assume $\dd \geq g-1$. Since $g=g_1+g_2+\dd -1$, this is equivalent to
$g_1+g_2\leq 2$.
There are thus four cases to consider:

1: $g_1=g_2=0;$

2: $g_1=0$, $g_2=1;$

3: $g_1=g_2=1$ and

4: $g_1=0$, $g_2=2.$

In the first three cases , i.e. when $g_i\leq 1$ for $i=1,2$, one checks that $(l,g-1-l)$ is indeed semistable for every
$l=0,\ldots,g-1$, so ${\mathcal T}=0$.

In  case 4  the multidegree $(g-1,0)$ is unstable, hence $\mathcal T$ is nontrivial
over the component $C_1^{g-1}\subset X^{g-1}\subset \X^{g-1}$ and hence $\agX$ is not natural.

Conversely, assume $\dd \leq g-2$.
We use examples~\ref{dcgvine} and \ref{Pvine}.
We have that $\dd $ is equal to the number of multidegree classes. Therefore the set
$\{(l, g-1-l),\  l=0,\ldots,g-1\}$, having cardinality $g\geq \dd+2$,  contains at least two  pairs of equivalent
multidegrees corresponding to two different classes.
Hence at least one of such pairs contains a nonsemistable multidegree, so the map $\ag$ is twisted
(i.e. ${\mathcal T}\neq 0$) and $\agX$ does depend on $f$.
\end{proof}
\end{nota}

\subsection{Abel maps of arbitrary degree}

Now we briefly illustrate what is known for a general $d\geq 1$.
\begin{nota}{\it Abel  maps of degree $1$.}
If $d=1$ the picture is rather well understood by \cite{AK},
\cite{EGK} for integral curves, and by \cite{CE} 
 for reducible ones.

With an approach  similar to what we used in case $d=g-1$, one constructs a
 map $\alpha^1_X:X\to \ov{P^1_X}$ which turns out to be natural and,  more remarkably,
regular. 
If $X$ is irreducible the choice of $\ov{P^1_X}$ is not an issue
(see \ref{N=Dirr}). 
 In case $X$ is reducible
 $\ov{P^1_X}$ can be chosen to be either the compactified Picard scheme constructed in
\cite{caporaso},
 or the one constructed in \cite{ed01}. We refer to the above mentioned  papers for details.

For a smooth curve $C$ (or   a family of smooth curves)
a fundamental fact is that the degree-$1$ Abel map is a closed embedding.
Moreover, recall that, after a point $p_0\in C$ is chosen, 
we can ``translate" $\alpha^1_C$ by composing it with the isomorphism
$\Pic^1C\to \Pic^0C$ mapping $[L]$ to $[L(-p_0)]$. In this way we get
the ``Abel-Iacobi" map
$$
\alpha_{p_0}:C\ha \Pic^0C=J(C)
$$ 
mapping $p\in C$ to $[\O(p-p_0)]$. Recall that $\alpha_{p_0}$ 
is a closed embedding, endowed with a universal property with respect
to mappings to Abelian varieties; namely
every map $h:C\to A$ of $C$ to an Abelian variety $A$
factors uniquely through $\alpha_{p_0}$, i.e. 
$h= h'\circ \alpha_{p_0}$ for a unique 
$h':A\to J(C)$.
In particular  $\alpha_{p_0}(C)$ generates $J(C)$
as a group.

What if $X$ is singular? Then the degree $1$-Abel map of $X$, and its  translates 
mapping to the compactified jacobian 
(see above),
 turn out to be  closed embeddings,
unless it cannot possibly be so for simple  reasons.
More precisely, it is not hard to see that if $X$ contains a smooth rational component $L$ that is attached to its
complement only in separating nodes of $X$, the degree-$1$ Abel map must contract   $L$  to a point
(note that $\Pic^dL$ is a point for all $d$).
These are the only cases where the map fails to be an embedding.
In particular, if $X$ is free from separating nodes, its degree-$1$ Abel map is a closed embedding.
\end{nota}

\begin{nota}{\it Abel maps of irreducible curves.}
\label{irrd}
The situation is more subtle for higher $d$, 
where not much is known about    compactifying Abel maps, even when there is no problem in defining them
simpy as rational map. For example,
 for irreducible curves there exists up to isomorphism a unique compactified Picard scheme and
one can proceed as we did for $d=g-1$. Now it is easy to see that
naturality is not an issue (see below).
So consider the rational map
$$
\alpha^d_X:X^d\da \PXb
$$
mapping $(p_1,\ldots,p_d )$ to $[\O_X(\sum p_i)]$ if all the $p_i$ are smooth points of $X$.
It is known that for $d\geq 2$ there is no hope that the map be   regular; indeed 
one needs to modify $X^d$ (blowing it up) to extend it. The explicit description of such a modification
is known in very  few cases,
for
$d=2$; see
\cite{coelho}.

A different approach which replaces $X^d$ with the Hilbert scheme of length $d$ subscheme on $X$ is pursued in
\cite{EK}.


\end{nota}

\begin{nota}{\it Naturality for arbitary $d$.} 
Let us now consider  reducible curves. 
In the previous section we have dealt  with the case $d=g-1$, 
considered Abel maps as rational maps,
and have seen that    there are very few
curves for which this map is natural  (see Proposition~\ref{g-1nat})

 For arbitrary  $d\geq
2$, the issue is complicated by the fact that, as we said, there exist different  compactified Picard 
schemes. Anyways,
after some choice for  $\PXb$  is made,
one  can work by specialization and produce  a rational map
\begin{equation}
\label{adX}
\adX: X^d \da \PXb
\end{equation}
similarly to what we did in case $g-1$.
We omit any explicit  definition to avoid choosing a specific compactified Picard scheme.
In fact, the point of this section is precisely   to describe  some facts that depend only on the type 
of   compactified Picard
scheme, and that can be proved without using the details of any specific  construction.
Namely, what about naturality of $\adX$?

We   shall concentrate on  compactified Picard schemes of N-type.
Observe that if $\pfb$ is of N\'eron type, the existence of a canonical map
$$
\ad:\X^d \da \pfb
$$
is guaranteed by the N\'eron mapping property. Indeed,
denote by $\dX^d\subset \X^d$ the smooth locus of $\X^d\to B$. 
The generic fiber $\gen^d$ of $\dX^d\to B$  has its own Abel map $\alpha_{\gen}^d:
\gen^d\to \Pic^d\gen$. The N\'eron mapping property yields a unique map 
$$
N(\alpha _{\gen}^d):\dX^d\la \nf.
$$
Now, as $\pfb$ is of N\'eron type, we can consider the map
$n_f^{-1}:\nf \to \pfb$ (cf. Definition~\ref{NDtype}). Composing, we get  a regular map
$$
\dX^d\  \stackrel{N(\alpha _{\gen}^d)}{\la} \  \nf\stackrel{n_f^{-1}}{\la} \pfb
$$
whose restriction to $X^d$ is the  rational map $\adX$ introduced in (\ref{adX}).
\end{nota}
We  need a definition:

\begin{defi}
\label{} Let $X$ be a (connected, nodal) curve and $\Gamma _X$ its dual graph.
The {\it essential  graph} of $X$ is the graph
$\ov{\Gamma _X}$ obtained from $\Gamma _X$ by eliminating every loop and by
contracting every separating edge to a point.

The edge connectivity of $\ov{\Gamma _X}$, (i.e.  the minimal number of edges that one needs to remove to
disconnect $\ov{\Gamma _X}$) will be called the {\it essential connectivity} of $X$ and denoted
$
\epsilon(X)$.
\end{defi}

\begin{example}
\label{ecvine}
If $X=\vine$ with $\dd\geq 2$ we have that $\Gamma_X=\ov{\Gamma _X}$;
if $\dd =1$ then $\ov{\Gamma _X}$ is a point. 
We have
\begin{displaymath}
\epsilon(\vine)=\left\{ \begin{array}{l}
\dd \  \text{ \  \  \  \  if } \dd\geq 2\\
+\infty \  \text{ if } \dd=1.\\
\end{array}\right.
\end{displaymath}
If $X$ is irreducible, $\ov{\Gamma _X}$ is a point and $\epsilon(X)=+\infty$.
\end{example}
The following  follows immediately from theorem 1.5 in \cite{CMM}.  
\begin{prop}
\label{natd}
Let $X$ be a curve and
 $\PXb$ be a compactified Picard scheme of N\'eron-type. Let
 $\alpha^d_{f,X}:X^d\da \PXb$ be the Abel map associated to a regular smoothing $f$ of $X$.
Then $\alpha^d_{f,X}$ is natural  only if $d<\epsilon (X)$.
\end{prop}

\begin{remark}
\label{large} Using Theorem~\ref{factMR} 
(which ensures that every curve admits some compactified Picard scheme on N-type)
one sees that the locus in $\mgbar$ of curves that fail to admit a natural Abel
map in degree $d\geq 2$ is quite large, i.e. it has codimension equal to 2.
Indeed, by Proposition~\ref{natd} and Example~\ref{ecvine},
the curve $\vine$ with $\dd = 2$ does not admit any natural $d$-th Abel map,
unless $d=1$.
\end{remark}

\begin{nota}{\it Abel maps for families over a higher dimensional base.}
Let $h:{\mathcal C}\to S$ be a family of smooth curves over any scheme $S$.
As we mentioned in \ref{smooth}, for every $d\geq 1$ there exists a (relative) $d$-th Abel map
$$
\alpha^d_h:{\mathcal C}^d\to \Pic^d_h=\Pic^d_{\mathcal C/S}
$$
where ${\mathcal C}^d$ denotes the fibered power over $S$ and $\Pic^d_h$ the relative Picard scheme in degree $d$
(a smooth projective scheme over $S$, as all fibers of $h$ are smooth).
In particular we can apply that to the universal family of smooth curves, and ask whether the construction extends 
to  stable curves.

More precisely, assume $g\geq 3$ and
let
$h_g:{\mathcal C}_g\to M_g^0$ (respectively $\ov{h_g}:\ov{{\mathcal C}_g}\to \mgbar^0$) be the universal family of
smooth (respectively stable) curves of genus $g$ over the  moduli space of automorphism-free curves.
Let ${\mathcal C}_g^d\to M_g^0$ and $\ov{{\mathcal C}_g}^d\to \mgbar^0$ be the $d$-th fibered powers
over the respective bases. Of course $\ov{h_g}$ is not a smooth morphism,
so, let us introduce its smooth locus,
denoted $\widetilde{h_g}:\widetilde{{\mathcal C}_g}\to \mgbar^0$ and its $d$-th fibered power
$$
\ov{{\mathcal C}_g}^d\supset \widetilde{{\mathcal C}_g}^d \la \mgbar^0.
$$
Consider now the $d$-th Abel map for the universal smooth curve
\begin{equation}
\label{unab}
\alpha^d_{h_g}:{\mathcal C}_g^d\la \Pic^d_{h_g}=\Pic^d_g;
\end{equation}
where $\Pic^d_g=\Pic_{{\mathcal C}_g/M_g^0}^d$ is a standard quick notation for the universal   Picard variety over
$M_g^0$. Assume now that $d$ is such that 
every $X\in \mgbar$ has a compactified degree-$d$ Picard scheme of N-type,
and that such compactified Picard schemes glue together over $\mgbar$; by Theorem~\ref{factMR}
this amounts to assume that $(d-g+1,2g-2)=1$.

Recall that  there exists a compactification $\pdbn\to \mgbar$ for $\Pic^d_g$, which is the moduli scheme for the stack
$\pdbstn$ (see \ref{factMR}).
So we ask: 
does the map (\ref{unab}) extend to a regular map 
$
\widetilde{{\mathcal C}_g}^d \la \pdbn
$?

From Proposition~\ref{natd}, one sees that, if $d\geq 2$, the answer is {\it no.}

By contrast, in case $d=1$  the answer is {\it yes}; in fact
 E. Esteves recently proved    that  the map
$\alpha^1_{h_g}$  extends  to a regular map 
over the whole of 
$\ov{\mathcal C}_g$;
and, more generally,  that this holds in its stack  version (\cite{estlett}).

Finally, a similar reasoning  works if $d=g-1$. In such a case we need Proposition~\ref{natg-1}
rather than Proposition~\ref{natd}, and the scheme $\overline{P_{g-1,g}}$ (see (\ref{g-1sch}))
parametrizing compactified Picard schemes in degree $g-1$. Again we obtain that 
{\it there exists no regular map
from $
\widetilde{{\mathcal C}_g}^{g-1}
$ to }$\overline{P_{g-1,g}}$.
\end{nota}
\subsection{The theta divisor of $\PXgb$.}
\begin{nota}{\it Theta divisor of a smooth curve.}
\label{thetasmooth} Let us consider   a smooth curve $C$ of genus $g\geq 2$. Using the set up
of  \ref{smooth}, the locus of effective line bundles in
 $\Pic^{g-1}C$ is a divisor
(by (\ref{AW})),   called the {\it theta divisor }  of $C$, and denoted
$$
\Theta(C):=W_{g-1}(C)\subset \Pic^{g-1}C.
$$
Many properties of the curve $C$ are encoded in the geometry of $\Theta(C)$.
For example, assume $g\geq 4$, then   $C$ is hyperelliptic iff
$\dim \Theta(C)_{\text{sing}} = g-3$ (where $\Theta(C)_{\text{sing}}$ is the singular locus of $\Theta(C)$);
if $C$ is not hyperelliptic, then $\dim \Theta(C)_{\text{sing}} = g-4$;
in both cases $\Theta(C)_{\text{sing}}$ has pure dimension and it is irreducible if $C$ is hyperelliptic.

On the other hand   $\Theta(C)_{\text{sing}}$   is precisely described in terms of
special line bundles on $C$ (a line bundle is called ``special" if its space of global sections has 
 dimension higher
than expected).
Indeed, the Riemann singularity theorem states that for every $L\in \Theta(C)$,
the multiplicity of $\Theta(C)$ at $L$ is equal to $h^0(C,L)$.
In particular we get 
$$
\Theta(C)_{\text{sing}} = W^1_{g-1}(C)=\{L\in \Pic^{g-1}C: h^0(C,L)\geq 2\}.
$$
Observe that, as we saw in (\ref{AW}), $\Theta(C)$ is the image of the $(g-1)$-th Abel map,
therefore it is irreducible.
Finally, recall  that $\Theta(C)$ is a principal polarization on $\Pic^{g-1}C$ (see \cite{ACGH}).
\end{nota}

\begin{nota}{\it Theta divisor of a generalized jacobian.}
Now suppose that $X$ is singular.
We may define, using the notation (\ref{Wmd}),
$$
\widetilde{\Theta(X)}:=\{L\in \Pic^{g-1}X: h^0(X,L)\geq 1\}=\coprod_{|\md|=g-1}\Wmd .
$$
We know already, by Proposition~\ref{path}, that   $\widetilde{\Theta(X)}$  
has only finitely many components of
the right dimension (i.e. of dimension $g-1$); indeed if $\md$ is not semistable we have $\Wmd = \Pic^{\md}X$. So $\widetilde{\Theta(X)}$ is not a divisor, if $X$ is reducible.

Assume that $\md$ is semistable, then 
 $\Wmd$ has dimension equal to $g-1$; moreover   the same holds for the image of the Abel map
$\amd$.
Since the Abel map is only a rational map, let us denote its closure by
$
\Amd :=\ov{\im  \amd}\subset \Pic^{\md}X.
$
We  now ask what the relation between $\Wmd$ and $\Amd$ is. Are they equal
 (as for a smooth
curve)?  If $\md$ is strictly semistable, the answer is {\it no}, as we shall see in the next example. 

\begin{example}
\label{ssred}
Let $X=\vine=C_1\cup C_2$ and consider the strictly semistable multidegree (notation in Example~\ref{vine})
 $$
\md=(g_1-1,g_2-1+\dd).$$
We already know that
$\Wmd$ has  dimension $g-1$ (by Prop.~\ref{path}).
We shall prove that $\Wmd$ has two different irreducible components, 
one of which (necessarily by Prop.~\ref{path}) coincides with $\Amd$.

Consider the normalization 
 $\nu:\Xn=C_1\coprod C_2\to X$. 
For any $L\in \Pic^{\md}X$ denote by $M=\nu^*L$. Since $\Xn$
is disconnected, $M$ is uniquely determined by its restrictions 
$L_1$ and $L_2$ 
to $C_1$ and $C_2$.
So, pick a pair $(L_1,L_2)\in \Pic^{(g_1-1,g_2-1+\dd)}\Xn$.
If $h^0(C_1,L_1)\neq 0$, i.e. if $L_1\in \Theta(C_1)$,
 then (as every  $L_2\in \Pic^{g_2-1+\dd}C_2$ has
$h^0(C_2,L_2)\geq\dd$ by Riemann-Roch) 
$$h^0(\Xn,M)=h^0(C_1,L_1)+h^0(C_2,L_2)\geq \dd+1.$$
  Therefore
for every $L\in \Pic X$ such that $\nu^*L=(L_1,L_2)=M$ we have that
$h^0(X,L)\geq h^0(\Xn, M)-\dd\geq 1$.
We conclude that
$\Wmd$ contains a closed subset $W_1$ given as
$$
W_1 = (\nu^*)^{-1}\bigr(\Theta(C_1)\times \Pic^{g_2-1+\dd}C_2\bigl),
$$
where
$$
\nu^*:\Pic^{\md}X\la \Pic^{g_1-1}C_1\times\Pic^{g_2-1+\dd}C_2
$$ 
is the pull-back map. The fibers of $\nu^*$ are irreducible of dimension
$\dd -1$ (by \ref{genjac}), hence $W_1$ is irreducible of dimension
$$
\dim W_1=(g_1-1)+g_2+(\dd -1)=g-1,
$$
so that $W_1$ is   an irreducible component of $\Wmd$.

Now suppose that $h^0(C_1,L_1)= 0$. Call $q_1,\ldots, q_{\dd}$ the points of $C_2$ mapping to the nodes
of $X$. If $h^0(C_2,L_2(-q_1-\ldots-q_{\dd}))\neq 0$,
every $L$ such that $\nu^*L=(L_1,L_2)$ lies in $\Wmd$.
The locus $D\subset \Pic^{g_2-1+\dd}C_2$
of line bundles $L_2 $ such that
$h^0(C_2,L_2(-q_1-\ldots-q_{\dd}))\neq 0$ 
is a so-called {\it translate} of $\Theta(C_2)$.
Indeed,
consider the isomorphism
 $u:\Pic^{g_2-1}C_2\la \Pic^{g_2-1+\dd}C_2$ mapping $N$ to $N(+q_1+\ldots+q_{\dd})$.
Then $D$ is equal to $u(\Theta(C_2))$. In particular $D$ is
 irreducible of dimension $g_2-1$.
We just constructed a second irreducible component $W_2$ of $\Wmd$:
$$
W_2=(\nu^*)^{-1}\bigr(\Pic^{g_1-1}C_1\times D \bigl).
$$
Arguing as for $W_1$, we see that $W_2$ is irreducible of dimension $g_1+g_2-1+\dd -1=g-1$.

We conclude that $\Wmd=W_1\cup W_2$ as wanted
(we leave it to the reader to show that in $\Wmd$ there is nothing else).

\end{example}
The previous example had $\md$ strictly semistable. 
If $\md$ is stable, 
$\Wmd$ turns out to be irreducible (see \cite{Ctheta})  and equal to
$\Amd$. 
\end{nota}

\begin{nota}{\it Compactifying the Theta divisor.}
The history follows the same pattern 
as for compactified Picard schemes and Abel maps, with
the case of irreducible curves being solved much earlier. As we mentioned (cf. Theorem~\ref{factg-1}),
the case of an irreducible curve was treated in \cite{soucaris} and \cite{esttheta},
where it is proved that, on the compactified jacobian of
$X$, there  exists a Cartier, ample divisor $\Theta(X)$  
 and that, just as for smooth curves, $3\theta(X)$ is very ample (in the latter paper).
For a reducible $X$, 
the situation remained open for  some time (despite some  breakthroughs in \cite{beau})
partly because of the diversity of the existing   compactified Picard schemes.
The proof of the fact that $\PXgb$ (the canonical one, see  \ref{g-1can}) has a 
theta divisor which is Cartier and ample is due to \cite{alex}. A new element in his construction 
is the use of
 semiabelian group actions, which also yields a placement of the pair 
($\PXgb$, $\Theta(X)$) within the degeneration theory of principally polarized abelian varieties
(see \cite{alex} section 5, in particular 5.4).

Comparing  to the smooth case, to the rich   picture we  partly sketched in \ref{thetasmooth}, the
subject   opens up now  with a variety of interesting issues and unsolved problems.
Indeed,  not much  is known about the geometry
of the theta divisor of a singular curve and about its interplay with the geometry of the curve and   its
jacobian.
\end{nota}
\begin{nota}{\it The Theta divisor of $\PXgb$.}
The definition of $\Theta(X)$ 
is given by the non vanishing of $h^0$. To do this properly we must say what the
points of $\PXgb$ parametrize. There are at least two good options: semistable
torsion free sheaves of rank $1$ and degree
$g-1$ (see \cite{alex} for example); or stable line bundles on the partial normalizations of $X$.
We will use the second option, for consistency with Proposition~\ref{Pstr},
and   because this description enables us to describe examples quite easily.

We stated in Proposition~\ref{Pstr} that a point in $\PXgb$ corresponds naturally to a line bundle on some
normalization of $X$, and we described  an instance of this correspondence in  \ref{strvine}.
We shall here use the same notation.
Let us denote by $\ell$ a point in $\PXgb$ and by $L\in \Pic^{\md}\XS$ the corresponding point
(recall that $\XS$ is the normalization of $X$ at $S\subset \sing$)
so that $\ell =\eSd(L)$.
We set
\begin{equation}
\label{h0}
h^0(\ell):=h^0(\XS, L)
\end{equation}

now we define 
\begin{equation}
\label{Theta}
\Theta(X):=\{\ell \in \PXgb: h^0(\ell)\neq 0\}.
\end{equation}
The fact  that the above set-theoretic definition coincides with the definition of the divisor $\Theta(X)$
studied in \cite{soucaris}, \cite{esttheta} and \cite{alex} is shown in \cite{Ctheta}. As we said, description
(\ref{Theta}) is good to easily give examplesas we are going to do next,
consistently with the stratification of Proposition~\ref{Pstr}.
 \end{nota}
\begin{example}
If $X$ is irreducible with $\delta$ nodes, and $g\geq 1$,
then $\Theta(X)$ is irreducible and
\begin{equation}
\label{Tirr}
\Theta(X)\cong W_d(X)\coprod 
\Bigr(\coprod_{i=1}^{\dd-1}\bigr(\coprod_{\stackrel{S_i\subset \sing}{ \#S_i=i}}W_{d-i}(\Xn_{S_i})\bigl)\Bigl)
\coprod\Theta(\Xn).
\end{equation}
Each stratum $W_{d-i}(\Xn_{S_i})$ is irreducible of dimension $g-1-i$.
\end{example}
\begin{example}
If $X=\vine$ with $\dd =1$ then (of course now $\Theta(X)$ is reducible)
\begin{equation}
\label{Tct}
\Theta(X_1)\cong \Theta(C_1)\times \Pic^{g_2-1}C_2\cup \Theta(C_2)\times \Pic^{g_1-1}C_1.
\end{equation}

If $X=\vine$ with $\dd =2$ then $\Theta(X)$ is irreducible and 
\begin{equation}
\label{Tban}
\Theta(X_2)\cong W_{(g_1,g_2)}(X) \coprod \Bigr( \Theta(C_1)\times \Pic^{g_2-1}C_2\cup \Theta(C_2)\times
\Pic^{g_1-1}C_1\Bigl).
\end{equation}
\end{example}

\

 \noindent
Lucia  Caporaso  caporaso@mat.uniroma3.it\  \\\
Dipartimento di Matematica, Universit\`a Roma Tre\\\ 
Largo S.\ L.\ Murialdo 1 \  \  00146 Roma - Italy

\end{document}